\documentclass[a4paper]{article}
\usepackage{amssymb}
\usepackage{authblk}
\usepackage{amsmath}
\usepackage{abstract}
\usepackage[utf8]{inputenc}
\usepackage[hidelinks]{hyperref}

\usepackage{titlesec}
\titleformat{\section}[block]{\center\scshape\large}{\thesection.}{0.7em}{}
\titleformat{\subsection}[block]{\center}{\thesubsection.}{0.7em}{}
\renewcommand\thesection{\Roman{section}}

\usepackage{amsthm}
\newtheorem{theorem}{Theorem}
\newtheorem{lemma}{Lemma}
\newtheorem{corollary}{Corollary}
\theoremstyle{definition}

\title{A Direct Proof of the Prime Number Theorem using Riemann's Prime-counting Function}

\author[a]{Zihao Liu}
\affil[a]{International Department, The Affiliated High School of SCNU,\authorcr Email: \url{mailto:travor_lzh@163.com}}
\date{}

\begin{document}
	\maketitle
	\begin{abstract}
		In this paper, we develop a novel analytic method to prove the prime number theorem in de la Vallée Poussin's form:
		$$
		\pi(x)=\operatorname{li}(x)+\mathcal O(xe^{-c\sqrt{\log x}})
		$$
		Instead of performing asymptotic expansion on Chebyshev functions as in conventional analytic methods, this new approach uses contour-integration method to analyze Riemann's prime counting function $J(x)$, which only differs from $\pi(x)$ by $\mathcal O(\sqrt x/\log x)$.
	\end{abstract}
	\section{Introduction}
	The prime number theorem\cite{apostol1976introduction}\cite{iwaniec_analytic_2004} has been a popular topic in analytic number theory since the 19th century. Its first proofs were independently given by Hadamard and de la Vallée Poussin in 1896 \cite{titchmarsh1986theory} using analytic methods. Since then, mathematicians such as Apostol \cite{apostol1976introduction}, Levinson \cite{levinson1966on}, Newman \cite{newman1980simple}\cite{zagier_newmans_1997}, Selberg \cite{selberg1949an}, Stein \cite{stein2003complex}, and Wright \cite{wright1952elementary} have explored different approaches (e.g. elementary methods, Tauberian theorems, and contour integrations) to prove this theorem. Instead of analyzing $\pi(x)$ directly, present proofs, regardless elementary or analytic, first attacked either of the Chebyshev functions
	\begin{equation}
		\label{eqche1}
		\vartheta(x)=\sum_{p\le x}\log p
	\end{equation}
	\begin{equation}
		\label{eqche2}
		\psi(x)=\sum_{n\le x}\Lambda(n)=\sum_{k\le\log_2x}\vartheta(x^{1/k})
	\end{equation}
	to derive asymptotic formulae, and estimates for $\pi(x)$ were later obtained using elementary methods, such as partial summation, to $\vartheta(x)$. Although Landau \cite{apostol1976introduction}\cite{landau1912uber} derived the prime number theorem from Möbius functions, he still used Chebyshev functions to build up the connection. Hence, it is an interesting task to investigate a method that avoids the use of \eqref{eqche1} and \eqref{eqche2}.

	In this paper, we propose a proof of the prime number theorem that avoids the use of Chebyshev functions. In particular, we first introduce Riemann's prime counting function $J(x)$ to our problem and relate it to $\pi(x)$. Then, by inverse Mellin transform, an integral representation is obtained for $\pi(x)$. In the subsequent part of this paper, we evaluate the contour integral and obtain the prime number theorem in the version of de la Vallée Poussin:
	\begin{theorem}[de la Vallée Poussin]
		\label{pnt}
		There exists a constant $c>0$ such that
		\begin{equation}
			\pi(x)=\operatorname{li}(x)+\mathcal O(xe^{-c\sqrt{\log x}})
		\end{equation}
		where $\operatorname{li}(x)$ is the logarithmic integral:
		\begin{equation}
			\label{eqli}
			\operatorname{li}(x)=\lim_{\varepsilon\to0^+}\left[\int_0^{1-\varepsilon}+\int_{1+\varepsilon}^x{\mathrm dt\over\log t}\right]
		\end{equation}
	\end{theorem}
	In essence, the main contribution of this paper is a new proof of the prime number theorem. The significance of this proof is that it is an analytic method that avoids the use of Chebyshev functions.
	\section{Riemann's prime counting function $J(x)$}
	Let $a_n$ be defined by
	\begin{equation}
		a_n=
		\begin{cases}
			1/k & n=p^k,p\text{ prime} \\
			0 & \text{otherwise}
		\end{cases}
	\end{equation}
	In his paper on $\pi(x)$, Riemann \cite{edwards1973riemann} defined $J(x)$ as the summatory function for $a_n$:
	\begin{equation}
		\label{eqjx}
		J(x)=\sum_{n\le x}a_n=\sum_{1<p^k\le x}\frac1k
	\end{equation}
	This function has a remarkable property that allows us to connect it to the standard prime counting function $\pi(x)$:
	\begin{theorem}
		\label{thmjxpi}
		For $x\ge2$, we have
		\begin{equation}
			\label{eqjxpi}
			J(x)=\pi(x)+\mathcal O(\sqrt x)
		\end{equation}
	\end{theorem}
	\begin{proof}
		It follows from \eqref{eqjx} that
		\begin{align*}
			J(x)-\pi(x)
			&=\sum_{2\le k\le\log_2x}{\pi(x^{1/k})\over k}
			\le{\sqrt x\over2}+\sum_{3\le k\le\log_2x}{\sqrt[3]x\over k} \\
			&\ll\sqrt x+\sqrt[3]x\log\log x
			\ll\sqrt x
		\end{align*}
	\end{proof}
	By \eqref{eqjxpi}, asymptotic formula for $J(x)$ can be ported to $\pi(x)$ with ease, so we can study $\pi(x)$ simultaneously when we analyze $J(x)$.
	\section{Integral representation for $J(x)$ and $\pi(x)$}
	Riemann \cite{edwards1973riemann} had shown that
	\begin{equation}
		\label{eqrmain}
		J(x)={1\over2\pi i}\int_{k-i\infty}^{k+i\infty}{x^s\over s}\log\zeta(s)\mathrm ds
	\end{equation}
	for noninteger $x>2$ and $k>1$ using Fourier analysis. However, because \eqref{eqrmain} involves infinite integrals, it would be difficult to analyze its asymptotic behavior directly, so we wish to convert \eqref{eqrmain} into a quantitative form. In other words, we prove a special case of Perron's formula \cite{titchmarsh1986theory}.
	\begin{lemma}
		\label{lemjz}
		For $\Re(s)>1$, we have
		\begin{equation}
			\label{eqjxzeta}
			\sum_{n=1}^\infty{a_n\over n^s}=\log\zeta(s)
		\end{equation}
		Moreover, we have
		\begin{equation}
			\label{eqnjdbound}
			\sum_{n=1}^\infty{|a_n|\over n^\sigma}\ll{1\over\sigma-1}
		\end{equation}
		as $\sigma\to1^+$.
	\end{lemma}
	\begin{proof}
		Since $\Re(s)>1$, the left hand side converges absolutely, we can change the order of summation safely to get
		\begin{align*}
			\sum_{n=1}^\infty{a_n\over n^s}
			&=\sum_p\sum_{k=1}^\infty{1\over kp^{ks}}
			=\sum_p\log{1\over1-p^{-s}} \\
			&=\log\prod_p{1\over1-p^{-s}}=\log\zeta(s)
		\end{align*}
		For \eqref{eqnjdbound}, we can use the fact that $|a_n|\le1$ and $\zeta(\sigma)\sim(\sigma-1)^{-1}$ as $\sigma\to1^+$ to justify it.
	\end{proof}
	\begin{theorem}
		For half-odd integer $x>2$ and $k=1+1/\log x$, we have
		\begin{equation}
			\label{eqpi}
			\pi(x)={1\over2\pi i}\int_{k-iT}^{k+iT}{x^s\over s}\log\zeta(s)\mathrm ds+R_1(x,T)
		\end{equation}
		with
		\begin{equation}
			R_1(x,T)\ll{x\log x\over T}+\sqrt x
		\end{equation}
	\end{theorem}
	\begin{proof}
		Plugging $|a_n|\le1$ and \autoref{lemjz} into Lemma 3.12 of \cite{titchmarsh1986theory}, we deduce
		\begin{equation}
			\label{eqjxr12}
			J(x)={1\over2\pi i}\int_{k-iT}^{k+iT}{x^s\over s}\log\zeta(s)\mathrm ds+\mathcal O\left(x\log x\over T\right)
		\end{equation}
		Finally, we replace $J(x)$ in \eqref{eqjxr12} by \eqref{eqjxpi} to obtain the result.
	\end{proof}
	\section{Evaluation of the integral}
	To handle the integral, we consider a rectangular region $R$ with vertices $1-\delta\pm iT$ and $k\pm iT$. This implies
	\begin{align}
		\label{eqi}
		{1\over2\pi i}\int_{k-iT}^{k+iT}{x^s\over s}\log\zeta(s)\mathrm ds
		&={1\over2\pi i}\oint_{\partial R}{x^s\over s}\log\zeta(s)\mathrm ds \\
		\label{eqi2}
		&+{1\over2\pi i}\left[\int_{k-iT}^{1-\delta-iT}+\int_{1-\delta-iT}^{1-\delta+iT}+\int_{1-\delta+iT}^{k+iT}\right]
	\end{align}
	where $\partial R$ denotes the boundary of $R$ in counterclockwise direction and $\delta$ is chosen such that $\zeta(s)\ne0$ in $R$ and on $\partial R$. This means that the only singularity of the integrand in $R$ is at $s=1$. However, because $s=1$ is a logarithmic singularity, circumvention is needed to evaluate it.
	\begin{lemma}
		If we define
		\begin{equation}
			\label{eqfr}
			f(r)={1\over2\pi i}\oint_{(r+)}{x^s\over s}\log{1\over s/r-1}\mathrm ds
		\end{equation}
		where $(r+)$ denotes any counterclockwise path that encloses only the singularity at $s=r$, then
		\begin{equation}
			f'(r)={x^r\over r}
		\end{equation}
		for $r\ne0$.
	\end{lemma}
	\begin{proof}
		Since
		\begin{align*}
			{\partial\over\partial r}\log{1\over s/r-1}
			={\partial\over\partial r}\left[\log r-\log(s-r)\right]
			=\frac1r+{1\over s-r}
		\end{align*}
		we have
		\begin{align*}
			f'(r)
			&={1\over2\pi i}\oint_{(r+)}{x^s\over s}\left[\frac1r+{1\over s-r}\right]\mathrm ds \\
			&=\lim_{s\to r}{x^s\over s}\left[{s-r\over r}+1\right]={x^r\over r}
		\end{align*}
	\end{proof}
	\begin{lemma}
		Using the same notation, for all $\eta>0$ we have
		\begin{align}
			\label{eqfr2}f(1+i\eta)=\operatorname{li}(x^{1+i\eta})-i\pi \\
			\label{eqfr3}f(1-i\eta)=\operatorname{li}(x^{1-i\eta})+i\pi
		\end{align}
	\end{lemma}
	\begin{proof}
		Because \eqref{eqfr} indicates that $f(r)$ vanishes when $\Re(r)\to-\infty$, we know that when $\eta>0$ there is
		\begin{align*}
			f(1+i\eta)
			&=f(1+i\eta)-f(-\infty+i\eta)=\int_{-\infty+i\eta}^{1+i\eta}f'(r)\mathrm dr \\
			&=\underbrace{\int_{-\infty+i\eta}^{1+i\eta}{x^r\over r}\mathrm dr}_{u=r\log x}
			=\int_{(-\infty+i\eta)\log x}^{(1+i\eta)\log x}{e^u\over u}\mathrm du
		\end{align*}
		Appealing to Cauchy's integral theorem \cite{titchmarsh2002theory}, we deform the path of integration of the last integral into a line segment connecting $-\infty$ and $-\varepsilon<0$, a clockwise circular arc connecting $-\varepsilon$ and $\varepsilon$, and finally a line segment connecting $\varepsilon$ and $(1+i\eta)\log x$. That is, we have
		\begin{equation}
			\label{eqfr4}
			f(1+i\eta)
			=\int_{-\infty}^{-\varepsilon}+\int_\varepsilon^{(1+i\eta)\log x}+\int_{-\varepsilon}^{\varepsilon}
		\end{equation}
		for all $\varepsilon>0$. As $\varepsilon\to0^+$, the latter integral becomes
		\begin{align*}
			\lim_{\varepsilon\to0^+}\int_{-\varepsilon}^{\varepsilon}{e^u\over u}\mathrm du
			&=\lim_{\varepsilon\to0^+}\int_\pi^0{e^{\varepsilon e^{i\theta}}\over \varepsilon e^{i\theta}}(i\varepsilon e^{i\theta})\mathrm d\theta=i\lim_{\varepsilon\to0^+}\int_\pi^0e^{\varepsilon e^{i\theta}}\mathrm d\theta \\
			&=-i\int_\pi^0\mathrm d\theta=-i\pi
		\end{align*}
		in which the interchanging of the limit operation follows from the fact that $e^{\varepsilon e^{i\theta}}$ converges to $1$ uniformly:
		\begin{equation}
			|e^{\varepsilon e^{i\theta}}-1|\le\sum_{n=1}^\infty{\varepsilon^n\over n!}=e^\varepsilon-1
		\end{equation}
		As $\varepsilon\to0^+$, the first two integrals become
		\begin{equation}
			\lim_{\varepsilon\to0^+}\underbrace{\left[\int_{-\infty}^{-\varepsilon}+\int_\varepsilon^{(1+\eta i)\log x}{e^u\over u}\mathrm du\right]}_{t=e^u}
			=\lim_{\varepsilon\to0^+}\left[\int_0^{e^{-\varepsilon}}+\int_{e^\varepsilon}^{x^{1+i\eta}}{\mathrm dt\over\log t}\right]
		\end{equation}
		By \eqref{eqli}, it is obvious that the right hand side evaluates to $\operatorname{li}(x^{1+i\eta})$. Plugging everything into \eqref{eqfr4}, we see that \eqref{eqfr2} is true. Adapting a symmetric argument yields \eqref{eqfr3}.
	\end{proof}
	\begin{theorem}
		\label{ti}
		Using the same notation, we have
		\begin{equation}
			{1\over2\pi i}\oint_{\partial R}{x^s\over s}\log\zeta(s)\mathrm ds=\operatorname{li}(x)
		\end{equation}
	\end{theorem}
	\begin{proof}
		Since $\zeta(s)\sim1/(s-1)$ when $s$ is near one \cite{titchmarsh1986theory}, we can transform the integrand:
		\begin{equation}
			{1\over2\pi i}\oint_{\partial R}{x^s\over s}\log\zeta(s)\mathrm ds
			={1\over2\pi i}\oint_{(1+)}{x^s\over s}\log{1\over s-1}\mathrm ds
		\end{equation}
		By \eqref{eqfr}, we see that the right hand side is exactly $f(1)$, which can be calculated by evaluating its Cauchy principal value from \eqref{eqfr2} and \eqref{eqfr3}:
		\begin{equation}
			f(1)=\lim_{\eta\to0^+}{f(1+i\eta)+f(1-i\eta)\over2}=\operatorname{li}(x)
		\end{equation}
	\end{proof}
	Since \autoref{ti} allows us to evaluate \eqref{eqi}, we now move our focus to \eqref{eqi2}. To begin with, we extract some properties of $\zeta(s)$ necessary for the proof.
	\begin{lemma}[de la Vallée Poussin]
		Let $s=\sigma+it$. Then there exists a fixed constant $A>0$ such that whenever $1-A/\log|t|\le\sigma\le2$ and $|t|\ge2$ we have $\zeta(s)\ne0$ and
		\begin{equation}
			\label{eqzeta}
			{\zeta'\over\zeta}(s)\ll\log|t|
		\end{equation}
		Moreover, in the same region we have
		\begin{equation}
			\label{eqlogzeta}
			|\log\zeta(s)|\ll\log|t|
		\end{equation}
	\end{lemma}
	\begin{corollary}
		\label{cologzeta}
		When $s\in\partial R$, we have
		\begin{equation}
			|\log\zeta(s)|\ll\log T
		\end{equation}
	\end{corollary}
	\begin{proof}
		Proof for \eqref{eqzeta} is already covered in the third chapter of \cite{titchmarsh1986theory}, so we move onto proving \eqref{eqlogzeta}. Since $|\log\zeta(2+it)|\le\log\zeta(2)$, we can rewrite $\log\zeta(s)$ into
		\begin{align*}
			\log\zeta(s)
			&=\int_2^\sigma{\zeta'\over\zeta}(u+it)\mathrm du+\mathcal O(1) \\
			&\ll\int_2^\sigma\log|t||\mathrm du|=|\sigma-2|\log|t|
		\end{align*}
		It is evident that $\sigma$ is bounded, so \eqref{eqlogzeta} directly follows.
	\end{proof}
	\begin{theorem}
		\label{ti2}
		If we choose $\delta=A/\log T$ and write
		\begin{equation}
			{1\over2\pi i}\int_{k-iT}^{k+iT}{x^s\over s}\log\zeta(s)\mathrm ds=\operatorname{li}(x)+R_2(x,T)
		\end{equation}
		then
		\begin{equation}
			R_2(x,T)\ll x\log^2T\exp\left(-{A\log x\over\log T}\right)
		\end{equation}
	\end{theorem}
	\begin{proof}
		By \autoref{cologzeta}, it is evident that the integral over the horizontal segments of $\partial R$ satsifies
		\begin{equation}
			\label{eqhbound}
			\left|\int_{k\pm iT}^{1-\delta\pm iT}{x^s\over s}\log\zeta(s)\mathrm ds\right|
			\ll{\log T\over T}\int_k^{1-\delta}x^u\mathrm du\ll{x\log T\over T\log x}
		\end{equation}
		Using \autoref{cologzeta}, we can also establish an upper bound for the left vertical segment of $\partial R$:
		\begin{align}
			\left|\int_{1-\delta-iT}^{1-\delta+iT}{x^s\over s}\log\zeta(s)\mathrm ds\right|
			&\ll x^{1-\delta}\log T\int_{-T}^T{\mathrm du\over|1-\delta+it|} \\
			\label{eqvbound}
			&\ll x^{1-\delta}\log^2T=x\log^2T\exp\left(-{A\log x\over\log T}\right)
		\end{align}
		wherein the last inequality follows from the definition of $\delta$. Plugging \eqref{eqhbound}, \eqref{eqvbound} into \eqref{eqi2} and applying \autoref{ti} to \eqref{eqi} give us the desired result.
	\end{proof}
	With everything prepared, we move our attention back to the prime counting function $\pi(x)$.
	\section{Proof of the \autoref{pnt}}
	Applying \autoref{ti2} to \eqref{eqpi}, we get
	\begin{equation}
		\label{eqpi2}
		\pi(x)=\operatorname{li}(x)+R_1(x,T)+R_2(x,T)
	\end{equation}
	If we set $\log T=\sqrt{\log x}$, then we have
	\begin{align}
		\label{eqr2}R_1(x,e^{\sqrt{\log x}})\ll x\log xe^{-\sqrt{\log x}} \\
		\label{eqr3}R_2(x,e^{\sqrt{\log x}})\ll x\log xe^{-A\sqrt{\log x}}
	\end{align}
	If we choose a proper $0<c<\min(1,A)$, then the logarithms in \eqref{eqr2} and \eqref{eqr3} can be absorbed into $e^{-c\sqrt{\log x}}$. Plugging them into \eqref{eqpi2} gives \autoref{pnt}.
	
	Using the fact that $\operatorname{li}(x)\sim x/\log x$, we can also sharpen the remainder of \eqref{eqjxpi}:
	\begin{corollary}
		Under \autoref{pnt}, we have for $x\ge2$ that
		\begin{equation}
			\pi(x)=J(x)+\mathcal O\left(\sqrt x\over\log x\right)
		\end{equation}
	\end{corollary}
	\begin{proof}
		Similar to how \autoref{thmjxpi} is proved, we just need to give a better upper bound for $J(x)-\pi(x)$. That is,
		\begin{align*}
			J(x)-\pi(x)
			&={\pi(x^{1/2})\over2}+\sum_{3\le k\le\log_2x}{\pi(x^{1/k})\over k} \\
			&\le{\pi(x^{1/2})\over2}+\pi(x^{1/3})\sum_{3\le k\le\log_2x}\frac1k \\
			&\ll{x^{1/2}\over\log x}+{x^{1/3}\log\log x\over\log x}\ll{\sqrt x\over\log x}
		\end{align*}
		where the second last $\ll$ follows from the fact that \autoref{pnt} implies $\pi(x)\ll x/\log x$.
	\end{proof}
	\section{Conclusion}
	In this paper, we present a new proof of the prime number theorem that avoids the use of Chebyshev functions. First, we prove \eqref{eqpi} to transform the arithmetic problem into an analytic problem. To evaluate the integral on the right hand side of \eqref{eqpi}, we set up a counterclockwise rectangular path so that this task is divided in to \eqref{eqi} and \eqref{eqi2}. Subsequently, differentiation under integral is applied to evaluate \eqref{eqi}, and classical zero-free region and upper bound for $\zeta(s)$ are used to give estimates for \eqref{eqi2}. Finally, choosing a proper $T$, we deduce again the prime number theorem in the version of de la Vallée Poussin.

	The most significant contribution of the proof is the evaluation of \eqref{eqi}, which gives the main term of \autoref{pnt}. Traditionally, number theorists study the prime number theorem and its generalizations (i.e. prime number theorem for arithmetic progressions\cite{montgomery2006multiplicative}, prime number theorem for automorphic L-functions\cite{liu_perrons_2007}, or sums of complex numbers over primes\cite{iwaniec_analytic_2004}) first using the weighted sum
	\begin{equation}
		\label{vsum}
		\sum_{n\le x}c_n\Lambda(n)=\sum_{p^k\le x}c_n\log p	
	\end{equation}
	and then convert it into sum over primes
	\begin{equation}
		\label{psum}
		\sum_{p\le x}c_p
	\end{equation}
	using partial summation\cite{apostol1976introduction}. However, introducing the weight $a_n$ suggests a possibility to directly estimate \eqref{psum} via
	\begin{equation}
		\label{ppsum}
		\sum_{n\le x}a_nc_n=\sum_{p^k\le x}{c_{p^k}\over k}=\sum_{p\le x}c_p+R
	\end{equation}
	To study \eqref{psum} from \eqref{vsum}, partial summation is needed for conversion, but studying \eqref{psum} from \eqref{ppsum} only requires us to provide upper estimates for
	\begin{equation}
		R=\sum_{2\le k\le\log_2x}\sum_{p\le x^{1/k}}{c_{p^k}\over k}
	\end{equation}
	Consequently, the proof technique presented in this paper may inspire new proofs of generalized prime number theorems\cite{liu_perrons_2007}\cite{iwaniec_analytic_2004} as they are also studied using Perron's formula\cite{titchmarsh1986theory}.
	\bibliographystyle{unsrt}
	\bibliography{refs}
	\newpage
	\section{Acknowledgment}
	Thanks Mr. Gui for supporting my project. I deeply appreciate my mother for providing me with a research environment. Particularly, she installed three white boards in my study room so that I can derive equations conveniently. Gratitude to all teachers and parents who have cared and aided me during this project!

\end{document}